\documentclass[11pt]{article}

\usepackage[centertags]{amsmath}
\usepackage{amsfonts}
\usepackage{amssymb}
\usepackage{theorem}
\usepackage{ulem}
\usepackage{isolatin1}
\usepackage{epsfig}
\usepackage{subfigure}
\usepackage[all]{xy}

\theoremstyle{break} \newtheorem{theorem}{Theorem}[section]
\theoremstyle{break} 
\theoremstyle{break}        
\theoremstyle{break} \newtheorem{lemma}[theorem]{Lemma}
\theoremstyle{break} 
\theoremstyle{break} \newtheorem{example}[theorem]{Example}
\theoremstyle{break} 
\theoremstyle{break}
{\theorembodyfont{\rmfamily}\newtheorem{remark}[theorem]{Remark}}
{\theorembodyfont{\rmfamily}}
\theoremstyle{break} 
\theoremstyle{break} 
\theoremstyle{break} 
\theoremstyle{break} 
\numberwithin{equation}{section}

\def\Re{\mathop{{\rm Re}}}

\begin{document}

\renewcommand{\thefootnote}{}%\fnsymbol{footnote}}
\stepcounter{footnote}
\begin{center}
{\bf \Large Critical points of 
inner functions, nonlinear partial differential\\[2mm] equations, and an extension of
Liouville's theorem
}
\end{center}
\renewcommand{\thefootnote}{\arabic{footnote}}
\setcounter{footnote}{0}
\begin{center}
{\large  \bf Daniela Kraus and Oliver Roth}\\[2mm]
{\small Universit\"at W\"urzburg,
 Mathematisches Institut, \\[-0.5mm]     
 D--97074 W\"urzburg, Germany}\\[1mm]
\end{center}

%\centerline{\today}

\bigskip
\renewcommand{\thefootnote}{\arabic{footnote}}

\smallskip
\begin{center}
\begin{minipage}{13cm} {\bf Abstract.}   
We establish an extension of Liouville's classical
representation theorem for solutions of the partial differential equation $\Delta u=4\, e^{2u}$ and combine this result
with methods from nonlinear elliptic PDE to
 construct holomorphic maps with prescribed critical points and specified
boundary behaviour. For instance, we show that for every Blaschke sequence
$\{ z_j\}$ in the unit disk there is always a Blaschke product with $\{z_j\}$
as its set of critical points. Our work is closely related to the
Berger--Nirenberg problem in differential geometry.
\end{minipage}
\end{center}

\section{Introduction}

In this paper we discuss a method for constructing holomorphic maps with
prescribed critical points based on a study of the Gaussian curvature equation
\begin{equation} \label{eq:gaussian0}
\Delta u=4 \, |h(z)|^2 \, e^{2 u} \, ,
\end{equation}
where $h$ is a holomorphic function on a domain $\Omega\subset \mathbb{C}$.
This technique has several applications to free  boundary value
problems for holomorphic maps of Riemann--Hilbert--Poincar\'e type.
For instance, we
prove the existence  of infinite Blaschke products with preassigned
  branch points satisfying the Blaschke condition. The construction is in two steps.
In a first step we find a solution of
the curvature equation (\ref{eq:gaussian0}) in the unit disk with  degenerate boundary data
$u=+\infty$ on the unit circle  when $h(z)$ is an infinite
Blaschke product (see Theorem  \ref{thm:pde} below). 
As we shall see, this is  a special case of the 
Berger--Nirenberg problem in differential geometry, i.e., the question which
functions $\kappa : R \to \mathbb{R}$ on a Riemann surface  $R$ arise as the
Gaussian curvature of a conformal Riemannian metric $\lambda(z) \, |dz|$ on
$R$. The Berger--Nirenberg problem is well--understood, when the Riemann surface 
$R$ is compact and not the sphere (see for instance   
Chang \cite{Chang}, and also  Moser 
\cite{Mos} and Struwe \cite{Str} for the case of the sphere), but is still 
not completely understood for  noncompact Riemann surfaces in which case
the Berger--Nirenberg problem is concerned with
{\it complete} conformal Riemannian metrics having prescribed curvature
 (see for instance Hulin \& Troyanov \cite{HT92} as one of the many
 references). 
In a second step, we establish an extension  
of Liouville's  classical  representation formula (see Liouville \cite{Lio1853}) for the solutions of
the Liouville equation
\begin{equation} \label{eq:liouville0}
\Delta u=4 \, e^{2 u}
\end{equation}
 to the more general equation 
(\ref{eq:gaussian0}) (see Theorem \ref{thm:liouville}). These two steps 
 combined with some standard results about bounded analytic functions 
allow a quick construction of the desired Blaschke product (Theorem \ref{thm:blaschke}).

\medskip

Building holomorphic maps with the help of the Liouville equation (\ref{eq:liouville0}) 
is an old idea and can be traced back at least to the work of Schwarz
\cite{Sch1891}, Poincar\'e \cite{Poi1898}, Picard
\cite{Pic1890,Pic1893,Pic1905} and Bieberbach \cite{Bie12,Bie16}. In fact, 
many of the  first attempts to prove the Uniformization Theorem for
Riemann surfaces were based on Schwarz' suggestion in
\cite{Sch1891}  to use  the partial differential  equation $\Delta u=4 \, e^{2 u}$ for
this purpose. After Poincar\'e and Koebe proved the Uniformization Theorem
by different means, the method seemed to have only  occasionally
been used  in
complex analysis. One notable important exception is M.~Heins celebrated paper
\cite{Hei62} in which the Schwarz--Picard problem (a special case of the
Berger--Nirenberg problem) was solved. We also note that a complete proof of the full
Uniformization Theorem via Liouville equation can be found in a recent
paper \cite{MT2002} by Mazzeo and Taylor. The main new aspect of the present
work is to show that the same method can also be applied  in situations when branch
points occur even though branching complicates the treatment considerably.
Roughly speaking, this is accomplished by replacing
 Liouville's equation (\ref{eq:liouville0}) by
 the Gaussian curvature equation (\ref{eq:gaussian0}) with the critical points  
encoded as the zeros of the holomorphic function $h(z)$.

\medskip

This paper is organized as follows. In Section \ref{sec:applications} 
we start with a discussion of some free and fixed boundary value problems for
analytic maps. 
Besides the construction of infinite Blaschke
products with preassigned critical points mentioned above, we also give the 
solution to a
problem raised by Fournier and Ruscheweyh \cite{FRold,FR01}  about
``hyperbolic'' finite Blaschke products, i.e., bounded analytic maps $f$
defined on a bounded simply connected domain $\Omega \subset \mathbb{C}$ such that
$$ \lim \limits_{z \to \partial \Omega} \frac{|f'(z)|}{1-|f(z)|^2}=1$$
with finitely many prescribed branch points.
These and the other results of Section \ref{sec:applications} are mainly
intended as examples 
illustrating the interplay between the Gaussian curvature
equation (\ref{eq:gaussian0}) and holomorphic functions with prescribed branching and
specified boundary behaviour -- the main topic of Section \ref{sec:method}.
There we discuss the basic ingredients we need for the proofs of the results of
Section \ref{sec:applications}: a solution of a special case of the
Berger--Nirenberg problem (Theorem \ref{thm:pde}) and an extension of
Liouville's theorem to the solutions of the variable curvature equation
(\ref{eq:gaussian0}) (see Theorem \ref{thm:liouville}). 
These results are proved in a final Section \ref{sec:proofs}, which also
includes a discussion of the necessary tools from nonlinear elliptic partial
differential equations and conformal geometry. 
In Section \ref{sec:method}
we also indicate how
methods from complex analysis can be used to some extent to obtain new
information about the Berger--Nirenberg problem. Thus the interaction between
the complex--analytic and the differential--geometric aspects of the curvature
equation works in both ways.
For instance, Theorem
\ref{thm:bergernirenberg} shows that there is always a  {\it unique} complete
conformal Riemannian metric on the unit disk  with curvature $\kappa(z)=-4 \,
|h(z)|^2$ when $h(z)$ is a Blaschke product.\footnote{On the other hand, 
when $h$ is a singular inner function,
  then in general there is more than one such metric, see Example \ref{ex:berger}.} 
Uniqueness results of this kind are usually
obtained by making use of Yau's generalized maximum principle \cite{Yau75,Yau78},
but require that the curvature is bounded above and below by negative
constants near the boundary (see Bland \& Kalka \cite{BK86} and Troyanov
\cite{Tro1992}). Theorem \ref{thm:bergernirenberg} allows instead infinitely
many zeros of the curvature function, which accumulate at the boundary. It
 hinges not only on  Yau's maximum principle, but also on 
a recent boundary version of
Ahlfors' lemma (see Kraus, Roth and Ruscheweyh \cite{KRR06}).

\section{Free and fixed boundary value problems for holomorphic maps with
  preassigned critical points} \label{sec:applications}

A proper holomorphic self--map of the unit disk $\mathbb{D}$ in the complex plane
$\mathbb{C}$ of degree $n$ has always a representation in terms of its $n$ zeros $a_1,
\ldots, a_n \in \mathbb{D}$ (with possible repetitions) as a 
 finite Blaschke product of the form
\begin{equation} \label{eq:blaschke}
 f(z)=\lambda \prod \limits_{j=1}^{n} \frac{z-a_j}{1-\overline{a_j} \,  z}\,
, \qquad |\lambda|=1 \, .
\end{equation}
 While for many questions in classical complex analysis
such a representation is perfectly well suited, it is not very helpful to 
study Blaschke products in their dependence on their critical
points. Knowledge about the
critical points of finite Blaschke products, however, is crucial in
a number of applications e.g.~when one studies the parameter space for complex polynomials 
(see Milnor \cite{Mil}). A second motivation for considering critical
points of Blaschke products comes from the new theory of discrete analytic
functions, where {\it discrete} Blaschke products are most naturally defined in
terms of their branch points (see Stephenson \cite{kS05}).

\medskip

 In fact finite Blaschke products are uniquely determined up to normalization by their
 critical points and these critical points can always be prescribed arbitrarily:

\medskip

{\bf Theorem A}

{\sl 
Let  $z_1, \ldots, z_n \in \mathbb{D}$ be given (not
necessarily distinct) points. Then there is always a
finite Blaschke product $f$ of degree $n+1$ with critical points $z_j$ and no
others. The Blaschke product $f$ is uniquely determined up to postcomposition
with a conformal automorphism of the unit disk.}

\medskip

The uniqueness statement in Theorem A is fairly straightforward and follows easily
for instance from Nehari's generalization of Schwarz' lemma
(see \cite[Corollary to Theorem 1]{Neh46}). 
To the best of our knowledge, the existence--part of
Theorem A was first proved by 
M.~Heins \cite[Theorem 29.1]{Hei62}. Heins' argument is
purely topological. He showed that
the set of critical points of all finite Blaschke products of degree $n+1$,
which is clearly closed, is also open in the polydisk $\mathbb{D}^n$ by using 
Brouwer's fixed point theorem.
Similar proofs were later given by
Q.~Wang \& J.~Peng \cite{WP79}, by T.~Bousch in his thesis \cite{B92}, and by 
S.~Zakeri \cite{Z96}. They consider the map $\Phi$ from (n+1)-tuples of
zeros to n-tuples of critical points (one degree of freedom
being used for normalization) and show that $\Phi$ is proper from the polydisk to the polydisk.
Again, with invariance of domain, this implies $\Phi$ is onto; see \cite{Z96} for
the details. Bousch shows that $\Phi$ is  even an analytic diffeomorphism \cite{B92}.
In particular, all these proofs are {\it nonconstructive}.

\medskip

An entirely novel and constructive approach to Theorem A based on Circle Packing
has recently been devised by Stephenson (see \cite[Theorem 21.1]{kS05}),  who builds
discrete finite Blaschke products with prescribed branch set and shows that they
 converge locally uniformly in $\mathbb{D}$ to a finite classical Blaschke product with specified
critical points. 

\medskip

The method we employ in the present paper for constructing
Blaschke products with prescribed critical points differs considerably from
the techniques described above and  is also  constructive in 
nature.\footnote{However, although our method is 
constructive, it is nevertheless not really suitable to 
   find a Blaschke product with prescribed critical points in an explicit
   form, even for finitely many critical points. Is there a 
   {\it finite} algorithm, which allows one to compute a finite Blaschke
   product from its critical points\,?} In addition, it 
has the advantage that it permits the construction of {\it infinite} Blaschke products
 with {\it infinitely} many prescribed critical points $z_1, z_2, \, \ldots$ 
provided $\{z_j\}\subset \mathbb{D}$ is a Blaschke sequence, i.e., it satisfies the Blaschke condition
\begin{equation} \label{eq:blaschkecondition}
 \sum \limits_{j=1}^{\infty} 1-\left| z_j\right|<\infty \, .
\end{equation} 
 Thus we have the following generalization of Theorem A.

\begin{theorem} \label{thm:blaschke}
Let $\{z_j\} \subseteq \mathbb{D}$ be a Blaschke sequence.
Then there exists a Blaschke product with critical points $\{z_j\}$
(counted with multiplicity) and no others.
\end{theorem}

Some remarks are in order. First note that unlike  Theorem A
there is no corresponding uniqueness statement in
Theorem \ref{thm:blaschke}. In fact, an infinite Blaschke product is not
necessarily determined by its critical points up to postcomposition
with a conformal disk automorphism. Here is a very simple example, when
there are no critical points at all.

\begin{example} \label{ex:frost}
Let $f$ be a Blaschke product which is also a
universal covering map of the unit disk onto a punctured disk $\mathbb{D}
\backslash \{ a \}$, $a\not=0$. For instance, one can take for $f$ any 
Frostman shift  $\tau_{\alpha} \circ F$
of the standard universal covering
$$ F : \mathbb{D} \to \mathbb{D} \backslash \{ 0 \} \, , \qquad F(z)=\exp \left( -
  \frac{1+z}{1-z} \right)$$
with a disk automorphism
$$ \tau_{\alpha}(z):=\frac{z+\alpha}{1+\overline{\alpha} \, z} \, , $$
provided $\alpha \in \mathbb{D} \backslash \{ 0 \}$. To check that $\tau_{\alpha}
\circ F$ is a Blaschke product for every $\alpha \in \mathbb{D} \backslash \{ 0 \}$ it
suffices to note that none of its angular limits is $0$, so it can have no
singular inner factor (see for instance \cite[Ch.~II, Theorem 6.2]{Ga}).
Now, $g(z):=f(-z)$ is also an infinite Blaschke product with the same critical points as
$f$, but clearly $f\not= T \circ g$ for any disk automorphism $T$.
%, since
%otherwise $f(0)=g(0)$ would imply $T=\text{id}$ and thus $f \equiv g$, which is certainly not true as $F(z)
%\not\equiv F(-z)$.
\end{example}

Secondly, for infinitely many branch points one cannot completely distinguish between
inner functions and Blaschke products in Theorem \ref{thm:blaschke}. 
In fact,  for
each inner function there are many Blaschke products with exactly the same 
critical points. This is an immediate consequence of  
 Frostman's theorem (\cite[Ch.~II, Theorem 6.4]{Ga}) that for any inner
function $F$ every Frostman shift $\tau_{\alpha} \circ F$ is a Blaschke
product for all $\alpha \in \mathbb{D}$ except for a set of capacity
zero.

\medskip

Thirdly, the condition in Theorem \ref{thm:blaschke} that $\{z_j\}$ is a Blaschke sequence
might be compared with a result of M.~Heins \cite{Hei62}, who showed that for 
any {\it bounded} nonconstant 
holomorphic map $f$ defined in $\mathbb{D}$ the critical points which are
contained in some fixed horocycle
$$ H(\omega,\lambda):=\left\{z \in \mathbb{D} \, : \left|1-z \, \overline{\omega}
\right|^2< \lambda \, \left(1-|z|^2 \right) \right\} \, , \qquad \omega \in
\partial \mathbb{D} \, , \, \lambda >0 \,  , $$
satisfy the Blaschke condition.\footnote{Even more is true: $f'$ restricted to any
horocycle is of bounded characteristic.} Hence there is a considerable gap
between the sufficient condition of Theorem \ref{thm:blaschke} 
for the critical points of an infinite Blaschke product 
and the above necessary condition of Heins. One is inclined to ask whether for
each nonconstant analytic map $f : \mathbb{D} \to \mathbb{D}$ there is always an inner
function or a Blaschke
product  with the same critical points as $f$.

\medskip

Fourthly, we wish to point out that the proof we give for Theorem
\ref{thm:blaschke} will show that in
the special case of finitely many branch points $z_1, \ldots , z_n \in \mathbb{D}$
there is always a {\it finite} Blaschke product (of degree $N$ say) with exactly
these critical points. The Riemann--Hurwitz formula then implies $N=n+1$. Thus
the existence--part of Theorem A might be considered as a special case 
of Theorem \ref{thm:blaschke}, see Remark \ref{rem:finite} below.

\medskip

The method we use to establish Theorem \ref{thm:blaschke}
is not restricted to the construction of Blaschke products. It can
also be used to prove the following extension of results due to
Fournier and Ruscheweyh \cite{FRold,FR01} and K\"uhnau \cite{Kuh} (see also
Agranovsky and Bandman \cite{AB96}).

\begin{theorem} \label{thm:4}
Let $\Omega \subseteq \mathbb{C}$ be a bounded simply connected domain, $z_1, \ldots, z_n$
finitely many points in $\Omega$, and
$\phi :\partial \Omega \to \mathbb{R}$  a continuous positive function. 
Then there exists a holomorphic function $f : \Omega \to \mathbb{D}$ 
with critical points at $z_j$ (counted with multiplicities) and no others such that
\begin{equation} \label{eq:bdd01}
 \lim \limits_{z \to \xi} \frac{|f'(z)|}{1-|f(z)|^2}=\phi(\xi)\, , \qquad \xi
\in \partial \Omega\, .
\end{equation}
If $g : \Omega \to \mathbb{D}$ is another holomorphic function with these properties, then
$g=T \circ f$ for some conformal disk automorphism $T : \mathbb{D} \to \mathbb{D}$. 
\end{theorem}

\begin{remark} \label{rem:rem}
\begin{itemize}
\item[(a)]
Choosing $\phi \equiv 1$, we see in particular that there is always a holomorphic solution 
$f : \Omega \to \mathbb{D}$ of the nonlinear
boundary value problem of {\it Riemann--Hilbert--Poincar\'e type}
\begin{equation} \label{eq:rusch}
\lim \limits_{ z \to \xi} \frac{|f'(z)|}{1-|f(z)|^2}=1 \, , \qquad \xi \in
\partial \Omega \, , 
\end{equation} 
with prescribed finitely many critical points in $\Omega$.
As indicated above, this
solves a problem which arises from the work in 
  \cite{FRold,FR01,Kuh}. There, Theorem \ref{thm:4} was proved
for the  special case $\phi \equiv 1$ and $\Omega=\mathbb{D}$  by using a completely different
method, which can be traced back to Beurling's celebrated extension of the
Riemann mapping theorem \cite{Beu53}. First, the problem is transfered to an 
integral equation which is then
solved by (i) iteration for a single critical point in \cite{Kuh} and by (ii)
applying Schauder's fixed point theorem  in the general case in \cite{FR01}.
In fact, an inspection of this method 
  shows that it requires at least some amount of regularity of $\Omega$\footnote{
For instance $\Omega$  of Smirnov--type (see \cite{Pom92}) would suffice.} 
and doesn't seem to be capable of yielding the full result
of Theorem \ref{thm:4} for general bounded simply connected domains.
\item[(b)] In general Theorem \ref{thm:4} does  not hold when $\Omega$ is not
  simply connected. See Example \ref{ex:ex} below.
\item[(c)]
Since (\ref{eq:bdd01}) and (\ref{eq:rusch}) are {\it free boundary value
  problems}  for the analytic map $f$ it 
is at first glance  a little surprising that one needs no assumptions on the
boundary regularity of $\Omega$ at all in Theorem \ref{thm:4}. 
The point is that  one can view  (\ref{eq:bdd01}) and (\ref{eq:rusch}) 
as a {\it fixed boundary value problem} for the conformal pseudo--metric
$$ \lambda(z) \, |dz|:=\frac{|f'(z)|}{1-|f(z)|^2} \, |dz| \, ,$$
i.e., $u(z):=\log \lambda(z)$ is a solution
 to Liouville's equation $\Delta u=4 e^{2 u}$ in $\Omega$ 
(except for the critical points) with {\it fixed} boundary values. 
In order to solve such a fixed boundary value problem it suffices that the
domain $\Omega$ is regular for the Laplace operator $\Delta$ in the sense of
potential theory.
\end{itemize}
\end{remark}

The boundary conditions (\ref{eq:bdd01}) and (\ref{eq:rusch}) involve  {\it
    unrestricted} approach to $\xi \in \partial \Omega$ from inside.
If $\Omega$ is a smooth domain, then we can
 relax this condition to nontangential limits and allow infinitely
many critical points for $f$. This is
the content of the following theorem, which we formulate for simplicity
 only for the
case $\Omega=\mathbb{D}$. Here and in the sequel we use the notation $\angle \lim$ to
indicate nontangential (angular) limits. 

\begin{theorem} \label{thm:5}
Let $\{z_j\}$ be a Blaschke sequence in $\mathbb{D}$
and let $\phi : \partial \mathbb{D} \to (0,\infty)$ be a function such that
$\log \phi \in L^{\infty}(\partial \mathbb{D})$.
Then there exists a holomorphic function $ f : \mathbb{D}\to \mathbb{D}$ with critical points
$z_j$ (counted with multiplicities)
such that
\begin{equation} \label{eq:bound}
 \sup \limits_{z \in \mathbb{D}} \frac{|f'(z)|}{1-|f(z)|^2} < \infty\, ,
\end{equation}
and 
$$ \angle \lim \limits_{z \to \xi} \frac{|f'(z)|}{1-|f(z)|^2}=\phi(\xi)
  \quad \text{ for a.e. } \xi \in \partial \mathbb{D}\, .$$
If $g : \mathbb{D} \to \mathbb{D}$ is another holomorphic function with these properties, then
$g=T \circ f$ for some conformal disk automorphism $T : \mathbb{D} \to \mathbb{D}$. 
\end{theorem}

\begin{remark} \label{prop:1}
Theorem \ref{thm:5} has the  following obvious partial converse.
If $f : \mathbb{D} \to \mathbb{D}$ is a non--constant holomorphic function such that
(\ref{eq:bound}) holds, then the non--tangential limit
$$  \angle \lim \limits_{z \to \xi}
\frac{|f'(z)|}{1-|f(z)|^2}=:\phi(\xi)$$
exists for a.e.~$\xi \in \partial \mathbb{D}$ and $\log \phi \in L^{1}(\partial
\mathbb{D})$. Moreover, the critical points of $f$ (counted with
multiplicity) satisfy the Blaschke condition.
To check this it suffices to say  that (\ref{eq:bound}) forces $f'$ to be
bounded on $\mathbb{D}$, so $f$ has a continuous extension to $\overline{\mathbb{D}}$.
\end{remark}

Theorems \ref{thm:blaschke}, \ref{thm:4} and \ref{thm:5} are all of a similar
flavour and will be proved in a unified way in Section \ref{sec:proofs}
below. However, there
are also a number of differences. For instance, unlike Theorem
\ref{thm:blaschke} we also have a uniqueness statement in
Theorem \ref{thm:5}. On the other hand, Theorems \ref{thm:blaschke} and
\ref{thm:5} deal with analytic maps defined on the unit disk, whereas
Theorem \ref{thm:4} is valid for any simply connected domain $\Omega$
regardless of the complexity of its boundary.

\section{The Gaussian curvature equation, critical points of holomorphic maps,
and the Berger--Nirenberg problem} \label{sec:method}

%\medskip

The idea of the proof of Theorem \ref{thm:blaschke} is based on the following
simple observation. If $f : \mathbb{D} \to \mathbb{D}$ is a bounded holomorphic map and the zeros
$\{z_j\}$ of $f'$ form a Blaschke sequence, then the Blaschke product
\begin{equation} \label{eq:blaschkeproduct}
 B(z):=\prod \limits_{j=1}^{\infty} \frac{-\overline{z}_j}{|z_j|} \frac{z-z_j}{1-\overline{z_j} \, z} \, 
\end{equation}
is a holomorphic self--map of $\mathbb{D}$. In particular,
\begin{equation} \label{eq:metric}
 \lambda(z) :=\frac{|f'(z)|}{1-|f(z)|^2} \frac{1}{|B(z)|} 
\end{equation}
defines the density of a conformal Riemannian metric $\lambda(z) \, |dz|$ on $\mathbb{D}$.
A quick computation shows that the Gaussian curvature
$$ \kappa_{\lambda}(z):=-\frac{\Delta \log \lambda(z)}{\lambda(z)^2}$$
of this metric is
$$ \kappa_{\lambda}(z)=-4 \, |B(z)|^2 \, .$$
Hence the metric $\lambda(z) \, |dz|$ has {\it nonpositive} curvature and the curvature
function $\kappa_{\lambda}$ is the negative square--modulus of a bounded holomorphic
function, which vanishes exactly at the critical points of $f$. 
In other words, the function
$$ u(z):=\log \lambda(z) \,  $$
is a smooth solution to the Gaussian curvature equation $\Delta u=4 \, |B(z)|^2 \, e^{2 u}$.
Note that the critical points $z_j$ are encoded as the zeros of the function
$B(z)$. Thus this PDE (or, what is the
same, the curvature of the metric $\lambda(z) \, |dz|$) is uniquely determined
by the critical points $\{z_j\}$ of $f$. 
If we assume momentarily that $f$ is a
finite Blaschke product, then $B$ is also a finite Blaschke product, and the
metric $\lambda(z) \, |dz|$ clearly blows up at the unit circle:
$$ \lim \limits_{z\to  \xi} \frac{|f'(z)|}{1-|f(z)|^2} \,
\frac{1}{|B(z)|}=+\infty \, , \qquad \xi \in \partial \mathbb{D} \, . $$
Thus the function $u(z)=\log \lambda(z)$
is a smooth solution to the boundary value problem
\begin{equation} \label{eq:bdd}
\begin{array}{rll}
 \Delta u & =4 \, |B(z)|^2 \, e^{2 u} & \text{ in } \mathbb{D} \, , \\ 
\lim \limits_{z \to \xi} u(z)& =+\infty &  \text{ for every } \xi \in \partial \mathbb{D}
\, .
\end{array}
\end{equation}

\medskip

The key idea of the proof of Theorem \ref{thm:blaschke} is now to reverse
these considerations.  Given a Blaschke sequence $\{z_j\}$ we form the Blaschke product
$B$ via (\ref{eq:blaschkeproduct}).
In a first step we then show that the boundary value problem
(\ref{eq:bdd}) always has a $C^2$--solution $u : \mathbb{D} \to \mathbb{R}$ 
(see Theorem \ref{thm:pde} below).
In a second step, we reconstruct from the corresponding conformal Riemannian
metric $\lambda(z) \, |dz|:=e^{u(z)} \, |dz|$ a holomorphic function
by solving the equation (\ref{eq:metric}) for $f$ (see Theorem
\ref{thm:liouville}).
In order to complete the proof of Theorem \ref{thm:blaschke} it 
finally remains to exploit the boundary condition in (\ref{eq:bdd}) 
to show that $f$ may be taken to be a Blaschke product.

\medskip

In fact, the two main steps of the proof of Theorem \ref{thm:blaschke} as
outlined above might both be stated in  more general form which we shall need
later. For instance, the boundary value problem (\ref{eq:bdd}) might be solved
for any bounded regular domain, that is, for any bounded 
domain possessing a Green's function which vanishes continuously on the boundary.

\begin{theorem} \label{thm:pde}
Let $\Omega \subseteq \mathbb{C}$ be a bounded  regular domain 
and $h : \Omega \to \mathbb{C}$, $h\not\equiv 0$,  a bounded holomorphic function. Then there exists a 
 $C^2$--solution $u : \Omega \to \mathbb{R}$ to 
\begin{equation} \label{eq:gaussian}
\Delta u = 4 \, |h(z)|^2 \, e^{2 u} \, \qquad \text{ in } \Omega 
\end{equation}
such that
\begin{equation}  \label{eq:boundary}
\lim \limits_{z \to \xi} u(z)=+\infty \text{ for every } \xi \in \partial \Omega
\, .
\end{equation}
\end{theorem}

\begin{remark} \label{rem:1}
Thus if $u : \Omega \to \mathbb{R}$ is a $C^2$--solution to the boundary value problem
(\ref{eq:gaussian})--(\ref{eq:boundary}), then the conformal Riemannian metric
$\lambda(z) \, |dz|:=e^{u(z)} \, |dz|$ has Gaussian curvature $-4 \, |h(z)|^2$.
If $\Omega$ is in addition a smooth domain (e.g.~$\partial \Omega$ is of class $C^2$), 
then the boundary condition (\ref{eq:boundary}) implies that this metric is even
a {\it complete} conformal Riemannian metric for $\Omega$. This follows e.g.~from the  boundary version
of Ahlfors' lemma in \cite{KRR06} (see in particular \cite[Theorem 5.1]{KRR06}).
Thus Theorem \ref{thm:pde} might be considered as a solution of the 
{\it Berger--Nirenberg problem} in the very special case that the curvature
of the metric $\lambda(z) \, |dz|$ is of the form $\kappa_{\lambda}(z)=-4 \,
|h(z)|^2$ for a bounded holomorphic function $h : \Omega \to \mathbb{C}$.
Related results have been obtained for instance by Bland \& Kalka \cite{BK86} and
Hulin \& Troyanov \cite{HT92}. They allow more general curvature
functions ($\kappa_{\lambda}(z)$ is only assumed to be H\"older continuous),
but they require $\kappa_{\lambda}(z)$ to be bounded below and above by negative
constants near the boundary. Theorem \ref{thm:pde} deals with a case of the Berger--Nirenberg
problem, when the curvature is only nonpositive, i.e., $\kappa_{\lambda}(z)$ is
allowed to vanish arbitrarily close to the boundary.
\end{remark}

Next, we consider the problem of uniqueness of  solutions   
 to the boundary value problem
(\ref{eq:gaussian})--(\ref{eq:boundary}). This  is a rather delicate matter.
In general the boundary value problem
(\ref{eq:gaussian})--(\ref{eq:boundary}) will have more than one solution, so
in view of Remark \ref{rem:1} 
there will be more than one complete conformal Riemannian metric on $\Omega$ having
curvature $-4 \, |h(z)|^2$. This is illustrated with the following example,
where the curvature comes from a singular inner function.

\begin{example} \label{ex:berger}
Let $h : \mathbb{D} \to \mathbb{D}$ be the singular inner function
$$ \exp \left(- \frac{1+z}{1-z} \right) \, . $$
Then a short calculation shows that
\begin{eqnarray*}
u_1(z) &=&\log \left[ \frac{1}{1-|z|^2} \, \frac{1}{|h(z)|} \right] \, \\[1mm]
u_2(z) &=& \log \left[ \frac{|h'(z)|}{1-|h(z)|^2} \frac{1}{|h(z)|} \right] 
\end{eqnarray*}
are two different solutions to (\ref{eq:gaussian})--(\ref{eq:boundary}), so
$e^{u_1(z)} \, |dz|$ and $e^{u_2(z)} \, |dz|$ are two complete conformal
Riemannian metrics on $\mathbb{D}$ with the same curvature.
\end{example}

At first glance, the nonuniqueness here might have been caused by the fact that the curvature $-4 \,
|h(z)|^2$ is not bounded from above by a {\it negative} constant. Indeed, 
Bland \& Kalka \cite{BK86}, Troyanov \cite{Tro1991}, and Hulin \& Troyanov 
\cite{HT92} have shown that if the curvature is bounded 
below and above by negative constants near the boundary of $\Omega$, then
there is at most one {\it complete} conformal Riemannian metric $\lambda(z) \, |dz|$ on
$\Omega$ with this curvature. These uniqueness results are essentially based on Yau's
generalized maximum principle \cite{Yau75,Yau78} which seems to require
an upper negative bound on the curvature at least near the boundary.
Our next result, however, gives a uniqueness result when the curvature is
allowed to vanish close to the boundary, even though it only deals with a
very special situation.

\begin{theorem} \label{thm:bergernirenberg}
Let $h$ be a Blaschke product. Then there exists a unique
complete regular\footnote{We call a conformal Riemannian metric $\lambda(z) \, |dz|$ on a domain $\Omega \subseteq
\mathbb{C}$ regular, if its density $\lambda : \Omega \to (0,+\infty)$ is of class
$C^2$.} conformal Riemannian metric $\lambda(z) \, |dz|$ on $\mathbb{D}$ with curvature
$-4 \, |h(z)|^2$. 
\end{theorem}

Theorem \ref{thm:bergernirenberg} and Example \ref{ex:berger} indicate that
even for the very special situation of the unit disk $\mathbb{D}$ and  curvature functions
of the form $-4 \, |h(z)|^2$ for  bounded holomorphic functions $h$ on $\mathbb{D}$, the 
question whether there is a unique complete conformal Riemannian metric 
with curvature $-4 \, |h(z)|^2$ seems to be quite intricate. We only mention
here that the uniqueness result of Theorem \ref{thm:bergernirenberg} can
easily be extended to bounded holomorphic functions $h$ without singular inner
factor in their canonical factorization (see \cite[Ch.~II,
Corollary 5.7]{Ga}), i.e.,
$$ h(z)=C \, B(z) \, S(z) \, , \qquad |C|=1 \, , $$
where $B$ is a Blaschke product, when we assume in addition that the outer
factor
$$ S(z)=\exp \left( \frac{1}{2 \pi}\int \limits_{0}^{2 \pi} \frac{e^{i
      \theta}+z}{e^{i \theta} -z} \, \log |f^*(e^{i \theta})| \, d \theta \right)$$
 is generated by an $L^1$--function $\log |f^*(e^{i \theta})|$ which is
 bounded away from $-\infty$.

\medskip

The above considerations all deal with particular cases of the Berger--Nirenberg
problem, when the curvature has the form $-4 \, |h(z)|^2$ for some holomorphic
function $h$. Theses cases appear to be very restrictive from 
 the viewpoint of differential geometry and nonlinear
partial differential equations. However, from the point of view of complex
analysis, they are most relevant for constructing holomorphic maps with
prescribed critical points.

\begin{theorem} \label{thm:liouville}
Let $\Omega \subseteq \mathbb{C}$ be a simply connected domain, $h : \Omega \to \mathbb{C}$ 
a
holomorphic function, $h \not\equiv 0$, and $u : \Omega \to \mathbb{R}$  a $C^2$--solution to
$\Delta u=4 \, |h(z)|^2 \, e^{2 u}$
in $\Omega$. Then there exists a holomorphic function $f : \Omega \to \mathbb{D}$ such
that
\begin{equation} \label{eq:liouville}
 u(z)=\log \left( \frac{|f'(z)|}{1-|f(z)|^2} \frac{1}{|h(z)|} \right) \, .
\end{equation}
Moreover, $f$ is uniquely determined up to postcomposition with a conformal
automorphism of the unit disk.
\end{theorem}

\begin{remark}
\begin{itemize}
\item[(i)]
For the special choice $h(z) \equiv 1$  Theorem \ref{thm:liouville} is a classical
result due to Liouville \cite{Lio1853}. The function $f$ is then locally
univalent and is sometimes called developing map of $u$.
\item[(ii)] Since $u$ is a $C^2$--solution, every holomorphic function
$f$ satisfying (\ref{eq:liouville}) has its critical points 
exactly at the zeros of the given
holomorphic function $h$ (counted with multiplicity).
\item[(iii)] There is a counterpart of Theorem \ref{thm:liouville} for
  conformal Riemannian metrics
  with {\it nonnegative} curvature: If $h$ is a holomorphic function on a
  simply connected domain $\Omega$, then every $C^2$--solution to $\Delta u=-
  \, 4 \, |h(z)|^2 \, e^{2 u}$ has a representation of the form
$$  u(z)=\log \left( \frac{|f'(z)|}{1+|f(z)|^2} \frac{1}{|h(z)|} \right) \, , $$
 where $f$ is a {\it meromorphic} function on $\Omega$. Moreover, $f$ is
 uniquely determined up to postcomposition with a rigid motion of the Riemann
 sphere. This result can be proved in exactly the same way as Theorem
 \ref{thm:liouville}.
\end{itemize}
\end{remark}

There are many proofs of Liouville's theorem  (the  case $h(z) \equiv 1$
of Theorem \ref{thm:liouville}) scattered throughout the literature (see for
instance Bieberbach \cite{Bie16}, Nitsche \cite{Nit57} and Minda \cite{Min})
and they continue to appear (see \cite{Bri}). They basically use the fact that
for $h(z) \equiv 1$ 
the Schwarzian derivative of the conformal Riemannian metric $e^{u(z)} \,
|dz|$, i.e.,
$$ \frac{\partial u^2}{\partial z^2}(z)-\left( \frac{\partial u}{\partial
    z}(z) \right)^2 \, , $$
 is a
{\it holomorphic} function $A(z)$ in $\Omega$ and $f$ can be found among the solutions
 to the Schwarzian differential equation 
$$ \left( \frac{f''(z)}{f'(z)} \right)'-\frac{1}{2} \left(
  \frac{f''(z)}{f'(z)} \right)^2=2 \, A(z) \, . $$
Since $A(z)$ is holomorphic in $\Omega$ {\it every} solution to this
differential equation is  meromorphic in $\Omega$. The main difficulty in
proving the more general Theorem \ref{thm:liouville} is that now $A(z)$
is meromorphic in $\Omega$ and indeed has poles of order $2$ exactly at the
zeros of $h(z)$. In the theory of Schwarzian differential equations (see Laine
\cite{Lai93}) this is known to be the most complicated situation. 
Nevertheless, it turns out that the associated Schwarzian differential equation still
has only meromorphic solutions, but this requires a considerable amount of work
({\it cf}.~Section \ref{sec:liouville}, in particular Lemma \ref{lem:liou2}).

\medskip

An immediate consequence of Theorem \ref{thm:liouville} is that the equation
$\Delta u=4 \, |h(z)|^2 \, e^{2 u}$ has no $C^2$--solution $u : \mathbb{C} \to \mathbb{R}$ if
$h(z)$ is an entire function. For the case  $h(z) \equiv 1$ this observation 
is due to Wittich \cite{Wit44}, Nitsche \cite{Nit57} and Warnecke \cite{War69}.

\section{Proofs} \label{sec:proofs}

\subsection{Proof of Theorem \ref{thm:blaschke}}
Let $\{z_j\}$ be a Blaschke sequence in $\mathbb{D}$. Then the corresponding Blaschke product
$$ B(z):=\prod \limits_{j=1}^{\infty} \frac{-\overline{z}_j}{|z_j|} \,
\frac{z-z_j}{1-\overline{z_j} \, z}$$
converges locally uniformly in $\mathbb{D}$, so $B : \mathbb{D} \to \mathbb{D}$ is a holomorphic
function. We employ Theorem \ref{thm:pde} with $h=B$ and get a $C^2$--solution
$u : \mathbb{D} \to \mathbb{R}$ of $\Delta u=4 \, |B(z)|^2 \, e^{2 u}$ in $\mathbb{D}$ with
$u=+\infty$ on $\partial \mathbb{D}$. Theorem \ref{thm:liouville} gives us a
holomorphic function $f : \mathbb{D} \to \mathbb{D}$ with 
$$ u(z)=\log \left( \frac{|f'(z)|}{1-|f(z)|^2} \frac{1}{|B(z)|} \right) \, .$$
Notice that the critical points of $f$ are exactly the zeros of $B$, i.e., the
points $z_j$, $j=1,2,\ldots$. We claim that $f$ is an inner function. This is
not difficult to see by considering the sets
\begin{eqnarray*}
A &:=& \left\{ \xi \in \partial \mathbb{D} \, \Big| \, f(\xi):=\angle \lim \limits_{z \to \xi}
  f(z) \text{ exists and } f(\xi) \in \partial \mathbb{D} \right\} \, , \\
A' &:= & \left\{ \xi \in \partial \mathbb{D} \, \Big| \, f(\xi):=\angle \lim \limits_{z \to \xi}
  f(z) \text{ exists and } f(\xi) \in  \mathbb{D} \right\} \, . 
\end{eqnarray*}
By Fatou's theorem \cite[Chapter I.5]{Ga} the union $A \cup A'$ has (one dimensional) Lebesgue measure $2
\pi$. If $\xi \in A'$, then
$$ \angle \lim \limits_{z \to \xi} \frac{|f'(z)|}{|B(z)|}=\angle  \lim
\limits_{z \to \xi}
\frac{|f'(z)|}{1-|f(z)|^2} \, \frac{1}{|B(z)|} \,  \left( 1-|f(z)|^2
\right)= \angle \lim \limits_{z \to \xi} \left[ e^{u(z)} \,  \left( 1-|f(z)|^2
\right) \right]=+\infty \,  $$
by construction. Thus the {\it holomorphic} function $B/f'$ has
angular limit $0$ at every point of the set $A'$. By Privalov's theorem
\cite[p.~94]{Ga} and
noting that $B/f'$ does not vanish identically, $A'$ is a nullset, so $A$ has
measure $2 \pi$ and $f$ is an inner function. Frostman's theorem
(see \cite[Ch.~II, Theorem 6.4]{Ga}) now guarantees that
$$ f_{\alpha}(z)=\frac{f(z)-\alpha}{1-\overline{\alpha} \, f(z)}$$
is a Blaschke product for all $\alpha \in \mathbb{D}$ except for a set of capacity
zero. In particular, $f_{\alpha}$ is a Blaschke product for some $|\alpha|<1$.
Since
$$ f_{\alpha}'(z)=\frac{1-|\alpha|^2}{\big(1-\overline{\alpha} \,
    f(z)\big)^2} \, f'(z) \, , $$
$f_{\alpha}$ is a Blaschke product with critical points $\{z_j\}$ and no
others. \hfill{$\square$}

\begin{remark} \label{rem:finite}
If $\{z_j\}$ is a finite sequence of $n$ points in $\mathbb{D}$, then $B(z)$ in the above proof is a finite
Blaschke product, so 
$$ \lim \limits_{z \to \zeta} \frac{|f'(z)|}{1-|f(z)|^2}=+\infty \qquad \text{ for
  every } \zeta \in \partial \mathbb{D} \,  $$
by construction.
Corollary 1.10 in \cite{KRR06} (which improves upon an earlier result of Heins
\cite{Hei86}) then implies that $f$ is a {\it finite} 
Blaschke product. The degree of this Blaschke product must be $n+1$ in view of
  the Riemann--Hurwitz formula and the existence--part of Theorem A follows.
\end{remark}

\subsection{Tools from PDE and conformal geometry} \label{sec:tools}

In this paragraph we collect some well--known facts from conformal geometry
and about the Gaussian curvature equation which we shall need in the sequel.

\medskip

The basic example of a conformal Riemannian metric 
is the {\it Poincar\'e  metric} $\lambda_{\mathbb{D}}(z) \, |dz|$ 
on the unit disk. Its density is given by
$$ \lambda_{\mathbb{D}}(z)=\frac{1}{1-|z|^2} \, . $$
When the domain $\Omega \subset \mathbb{C}$ has at least two boundary points, usually
dubbed hyperbolic domain, and $-a$ is a negative constant, then $\Omega$ carries a
unique complete regular conformal Riemannian metric with constant curvature $-a$. 
This metric  $\lambda(z) \, |dz|$
is obtained from the Poincar\'e metric $\lambda_{\mathbb{D}}(z) \, |dz|$ by means of
a universal cover projection $\pi : \mathbb{D} \to \Omega$ from
$$ \lambda(\pi(z)) \, |\pi'(z)| =\frac{2}{\sqrt{a}} \, \lambda_{\mathbb{D}}(z)=
\frac{2}{\sqrt{a}} \frac{1}{1-|z|^2} \, .$$
We call $\lambda(z) \, |dz|$ the hyperbolic metric on $\Omega$ with curvature $-a$.
Unless explicitly stated otherwise we usually take the normalization $a=4$ and
call the corresponding metric {\it the} hyperbolic metric of $\Omega$ (with
constant curvature $-4$). This metric is denoted by $\lambda_{\Omega}(z) \, |dz|$.

\medskip

An important result about conformal Riemannian metrics is the   {\it Ahlfors
  lemma} \cite{Ahl38}. It states that
the hyperbolic metric $\lambda_{\Omega}(z) \, |dz|$ is maximal in the sense
that for every regular conformal Riemannian metric $\lambda(z) \, |dz|$ with curvature
$\le -4$ the inequality $\lambda(z) \le \lambda_{\Omega}(z)$ holds for every
point $z \in \Omega$. This is even true for regular conformal {\it pseudo--metrics}
$\lambda(z) \, |dz|$, i.e., the density $\lambda$ is not necessarily strictly
positive, but is only assumed to be nonnegative, and $\lambda$ is of class
$C^2$ off its zero set.
On the other hand, for every {\it complete} regular conformal Riemannian {\it
  metric} on $\Omega$
with curvature bounded below by $-4$ the estimate $\lambda(w) \ge
\lambda_{\Omega}(w)$  is valid in $\Omega$. We refer to this fact as the
Ahlfors--Yau lemma (see Yau \cite{Yau75,Yau78}). There are boundary versions
of these results ({\it cf}.~Bland \cite{Bla83}, Troyanov \cite{Tro1991}, and Kraus,
Roth \& Ruscheweyh \cite{KRR06}).

\medskip

We next compile a number of facts about the Gaussian curvature equation
(\ref{eq:gaussian}). We adopt standard notation, so $C(\overline{\Omega})$ is the
set of real--valued continuous functions on the set $\overline{\Omega} \subset
\mathbb{C}$ and $C^k(\Omega)$ is the set of real--valued functions having all derivatives
of order $\le k$ continuous in the open set $\Omega \subseteq \mathbb{C}$.

\begin{lemma} \label{lem:help}
Let $\Omega \subseteq \mathbb{C}$ be a bounded  regular domain and $h : \Omega \to \mathbb{C}$  a
bounded holomorphic function. 
\begin{itemize}
\item[(a)] (Comparison principle)

If $u_1, u_2 \in C^2(\Omega) \cap C(\overline{\Omega})$  are two solutions
to (\ref{eq:gaussian}) and if $u_1 \le u_2$ on  the boundary $\partial \Omega$, then
$u_1 \le u_2$ in $\overline{\Omega}$. 
\item[(b)] Let $g_{\Omega}(z,\zeta) \ge 0$ denote Green's function of the regular domain
$\Omega$.  If $u$ is a bounded and integrable  function on  $\Omega$, then
$$v(z):=  -\frac{1}{2 \pi} \iint  \limits_{\Omega} g_{\Omega}(z, \zeta) \, 
\, 4 \, |h(\zeta)|^2 \, e^{2 u(\zeta)} \, d\!\; \sigma_{\zeta}, \qquad z \in \Omega\,, $$
where $\sigma_{\zeta}$ denotes two--dimensional Lebesgue measure,
belongs to $C^{1}(\Omega) \cap  C(\overline{\Omega})$ 
and  $v \equiv 0$ on $\partial \Omega$. If, in addition, $u$ is locally H\"older continuous with
exponent $\beta$, $0< \beta \le 1$, then $ v \in C^2(\Omega)$ and 
 $\Delta v=4 \, |h(z)|^2 \, e^{2u}$ in $ \Omega$.
\item[(c)]
If $u : \Omega \to \mathbb{R}$ is a $C^2$--solution to (\ref{eq:gaussian}) 
 and  $u$ is continuous on the closure
  $\overline{\Omega}$, then 
the integral formula
\begin{equation} \label{eq:green}
 u(z)=H(z)-\frac{1}{2 \pi} \iint  \limits_{\Omega} g_{\Omega}(z, \zeta) \, 
\, 4 \, |h(\zeta)|^2 \, e^{2 u(\zeta)} \, d\!\; \sigma_{\zeta}\, 
\end{equation}
holds for every $z \in \Omega$. Here, $H$ is harmonic in $\Omega$ and
continuous on $\overline{\Omega}$ with boundary values $u$, i.e.,
$$H \big|_{\partial \Omega} \equiv u \big|_{\partial \Omega}\, .$$
% and $d\!\;
%\sigma_{\zeta}$ denotes two--dimensional Lebesgue measure.
Conversely, if $u$ is a locally integrable and bounded function on the regular
domain $\Omega$  satisfying (\ref{eq:green}) for some
harmonic function $H$ in $\Omega$ which is continuous on $\overline{\Omega}$, then
$u$ belongs to $C^2(\Omega) \cap C(\overline{\Omega})$ and is a solution to
$\Delta u=4 \, |h(z)|^2 \, e^{2u}$ in $\Omega$ with  $u  \equiv H$ on $\partial
\Omega$.
\item[(d)] If $\phi : \partial \Omega \to \mathbb{R}$ is a continuous function, then 
there exists a
unique $C^2$--solution $u : \Omega \to \mathbb{R}$ continuous on $\overline{\Omega}$
to  the boundary value problem
\begin{alignat*}{2}
\Delta u& =4 \, |h(z)|^2 \, e^{ 2u}  \qquad& \text {in}\, \, \,  &\Omega \\
       u&=\phi                     \quad &\text{on} \, \, \, &\partial \Omega\, .
\end{alignat*}
\end{itemize}
\end{lemma}

For the comparison principle (Lemma \ref{lem:help} (a)) see \cite[Theorem
10.1]{GT97}. Part (b) of Lemma \ref{lem:help} is established 
in \cite[p.~54/55]{GT97} and \cite[p.~241]{CH68}, and part (c) follows by
combining (a) and (b). Lemma \ref{lem:help} (d) finally might be found in
\cite{GT97}, in particular Theorem 12.5 and the remarks on p.~308/309.

%\bigskip

We now can also give an example which illustrates that Theorem \ref{thm:4}
does not remain valid in general, when the domain $\Omega$ is not simply connected,
see Remark \ref{rem:rem} (b).

\begin{example} \label{ex:ex}
Let $\Omega$ be the annulus $\{z \in \mathbb{C} \, : \, 1/4
  <|z|<1/2\}$ and let $\phi : \partial \Omega \to \mathbb{R}$ be the continuous function
$$ \phi(\xi)=\left\{ \begin{array}{c}  \displaystyle \frac{4}{3} \\[4mm]
\sqrt{2} \end{array} \quad \text{ if } \quad \begin{array}{l} 
 |\xi|=\displaystyle \frac{1}{4} \, , \\[3mm]
 |\xi|=\displaystyle \frac{1}{2} \, .
\end{array} \right. $$
Now assume there exists a locally univalent holomorphic
function $f : \Omega \to \mathbb{D}$ such that 
$$ \lim \limits_{z \to \xi} \frac{|f'(z)|}{1-
|f(z)|^2} =\phi(\xi) \qquad \text{ for every } \, \,  \xi
\in \partial \Omega\, .$$ 
Note,
$$ \tilde{\lambda}(z)=\frac{1}{2 \sqrt{|z|} \, (1-|z|)}$$
is the density of a conformal Riemannian metric  of constant curvature $-4$ in $\Omega$ and
$\tilde{\lambda}(\xi)=\phi(\xi)$ for $\xi \in \partial \Omega$.
Hence, by the uniqueness part of Lemma \ref{lem:help} (d), we see that
$$ \tilde{\lambda}(z)=\frac{|f'(z)|}{1-|f(z)|^2}, \qquad z \in \Omega.$$
On the other hand, in the simply connected domain $D=\Omega \backslash (-1/4,-1/2)$, 
we have $$\tilde{\lambda}(z)=\frac{|g'(z)|}{1-|g(z)|^2}$$ for the holomorphic
function $g : D \to \mathbb{D}$, $g(z)=\sqrt{z}$.
Applying Theorem \ref{thm:liouville}, we deduce $g=T \circ f$ in $D$ for
some automorphism $T$ of $\mathbb{D}$. Thus $g$ has an analytic extension to
$\Omega$ given by $T \circ f$ which is absurd. 
\end{example}

\subsection{Proof of Theorem \ref{thm:pde}  and Theorem
  \ref{thm:bergernirenberg}}

{\bf Proof of Theorem \ref{thm:pde}.}
In order to prove Theorem \ref{thm:pde}, we consider for each integer $n \ge 1$ the
unique real--valued solution $u_n \in C^2(\Omega) \cap C(\overline{\Omega})$ 
to the boundary value problem 
\begin{alignat*}{2}
\Delta u& =4 \, |h(z)|^2 \, e^{ 2u}  \qquad& \text {in}\, \, \,  &\Omega \\
       u&=n                     \quad &\text{on} \, \, \, &\partial \Omega\, ,
\end{alignat*}
see Lemma \ref{lem:help} (d). By the comparison principle (Lemma
\ref{lem:help} (a)) we get a monotonically increasing sequence $\{ u_n\}$.
Our task is to show that $u_n$ converges to a solution of the boundary value
problem (\ref{eq:gaussian})--(\ref{eq:boundary}).
We proceed in a series of steps.

\begin{itemize}
\item[(i)] For each fixed $n \ge 1$ consider the conformal pseudo--metric
$$ \lambda_n(z) \, |dz|:= |h(z)| \, e^{u_n(z)} \, |dz| \, . $$ 
The Gaussian curvature of $\lambda_n(z) \, |dz|$ is found to be
$$ \kappa_{\lambda_n}(z)=- \frac{\Delta u_n(z)}{|h(z)|^2 \, e^{2
    u_n(z)}}=-4 \, . $$
Thus, by Ahlfors' lemma,
\begin{equation} \label{eq:ahlfors}
 \lambda_n(z) \le \lambda_{\Omega}(z) \, , \qquad z \in \Omega \, ,
\end{equation}
so we get a uniform bound
\begin{equation} \label{eq:boundabove}
u_n(z) \le \log \lambda_{\Omega}(z)-\log |h(z)| \, , \qquad z \in
\Omega \, .
\end{equation}
In particular, $\{u_n\}$ converges monotonically to some limit function $u :
\Omega \to \mathbb{R} \cup \{ +\infty\}$. Note that $u(z)$ is certainly finite when $z$ is not a
zero of $h$.
We need to show that $u \in C^2(\Omega)$, $\Delta u=4 \,
|h(z)|^2 \, e^{2 u}$ in $\Omega$ and $u(z) \to +\infty$ whenever $z \to \xi
\in \partial \Omega$.
\item[(ii)] Let $G$ be a regular domain, which is compactly contained in
  $\Omega$. Then by Lemma \ref{lem:help} (c),
\begin{equation} \label{eq:integral}
 u_n(z)=H_n(z)-\frac{1}{2 \pi} \iint  \limits_{G} g_{G}(z, \zeta) \, 
\, 4 \, |h(\zeta)|^2 \, e^{2 u_n(\zeta)} \, d\!\; \sigma_{\zeta}\,
\end{equation}
with a harmonic function $H_n : G \to \mathbb{R}$ which is continuous on
$\overline{G}$ and $H_n \equiv u_n$ on $\partial G$.
Note that the second term on the right--hand side is uniformly bounded in
$G$. This follows from
\begin{eqnarray*}
 0 & \le &  \frac{1}{2 \pi} \iint  \limits_{G} g_{G}(z, \zeta) \, 
\, 4 \, |h(\zeta)|^2 \, e^{2 u_n(\zeta)} \, d\!\; \sigma_{\zeta} = 
\frac{1}{2 \pi} \iint  \limits_{G} g_{G}(z, \zeta) \, 4  
\, \lambda_n(\zeta)^2  \, d\!\; \sigma_{\zeta} \\
&\overset{(\ref{eq:ahlfors})}{\le}& \frac{1}{2 \pi} \iint  \limits_{G} g_{G}(z, \zeta) \, 
\, 4 \, \lambda_{\Omega}(\zeta)^2  \, d\!\; \sigma_{\zeta} \le C_1(G) \cdot \frac{1}{2 \pi} \iint  \limits_{G} g_{G}(z, \zeta) \, 
\,   \, d\!\; \sigma_{\zeta} \le C(G) <\infty\, ,   
\end{eqnarray*}
where we have used 
the facts that the hyperbolic density $\lambda_{\Omega}$ is clearly
bounded on $G$ and 
$$ v(z):=  \frac{1}{2 \pi} \iint  \limits_{G} g_{G}(z, \zeta) \, 
\,   \, d\!\; \sigma_{\zeta}$$
is the solution to the boundary value problem $\Delta v=-1$ in $G$ and $v=0$
on $\partial G$, so $v \in C(\overline{G})$ is in particular bounded 
on $G$.
\item[(iii)] We next show that the sequence $\{H_n\}$ of harmonic functions in
  $G$ converges locally uniformly in $G$. Since $H_n\equiv u_n$ on $\partial
  G$ and $\{u_n\}$ is monotonically increasing on $G$,
it suffices by Harnack's theorem to show
  that $\{H_n(z_0)\}$ is bounded above at some point $z_0 \in G$. 
But this follows immediately from equation (\ref{eq:integral}),
 by what we have proved in part (ii)  and noting that (\ref{eq:boundabove})
implies $\{u_n(z_0)\}$ is
 bounded above if $h(z_0) \not=0$.

\item[(iv)] From (iii) we deduce that $\{H_n\}$ converges locally uniformly in
  $G$ to some harmonic function $H: G \to \mathbb{R}$. Since we have proved in part (i)
 that $\{u_n\}$ converges monotonically in $\Omega$ to some limit function
$u$, we thus see from part (ii) and formula (\ref{eq:integral})  that
$u(z)$ is finite for every $z \in G$ and
 $$  u(z)=H(z)-\frac{1}{2 \pi} \iint  \limits_{G} g_{G}(z, \zeta) \, 
\, 4 \, |h(\zeta)|^2 \, e^{2 u(\zeta)} \, d\!\; \sigma_{\zeta}\, , \qquad z
\in G \, $$
by Lebesgue's theorem on monotone convergence.
Now Lemma \ref{lem:help} (c) shows that $u$ belongs to $C^2(G) \cap C(\overline{G})$
and solves (\ref{eq:gaussian}) in $G$.  Since
$G$ is an arbitrary regular domain compactly contained in $\Omega$, we see that $u \in
C^2(\Omega)$ solves (\ref{eq:gaussian}) in $\Omega$.
\item[(v)] It is  now very easy to verify the boundary condition 
(\ref{eq:boundary}). Assume to the contrary 
that we can find a sequence $\{z_j\}\subset \Omega$ which converges to
$\xi \in \partial \Omega$ and a constant $0<C_2<\infty$ such that $u(z_j) <C_2$
for all $j$. Now choose an integer $m> C_2$.
Since $u_m(z_j) \to m$ as $j \to \infty$, there is an integer
 $J$ such that $u_m(z_j)>C_2$
for all $j>J$. But then the monotonicity of $\{u_n\}$ yields
$u(z_j) \ge C_2$ for all $j>J$ and the contradiction is apparent. Thus
$$\lim \limits_{z \to \xi} u(z)=+\infty$$
for every $\xi \in \partial \Omega$ as desired. \hfill{$\square$}
\end{itemize}
\hfill{$\blacksquare$}

{\bf Proof of Theorem \ref{thm:bergernirenberg}.}
The existence--part of Theorem \ref{thm:bergernirenberg} follows from Theorem 
\ref{thm:pde} and Remark \ref{rem:1}, so we need only show the uniqueness statement.
Let $\lambda(z) \, |dz|$ be a complete regular conformal metric on
$\mathbb{D}$ with curvature $-4 \, |h(z)|^2$, where $h$ is a Blaschke product, so its
curvature is  bounded from below by $-4$. The Ahlfors--Yau lemma
 then implies that
its density $\lambda(z)$ is bounded below by $1/(1-|z|^2)$, the 
 density of the hyperbolic metric on $\mathbb{D}$
with constant curvature $-4$. Thus $u(z):=\log \lambda(z)$ is a $C^2$--solution
to the boundary value problem (\ref{eq:gaussian})--(\ref{eq:boundary}).
Let $u_1$ and $u_2$ be two such solutions. We need to show $u_1 \equiv u_2$.
Consider the auxiliary function
$$ v(z):=\max\{u_1(z),u_2(z)\}-u_2(z) \, . $$
We first notice that $v(z)$ is subharmonic in $\mathbb{D}$. In fact,
if $v(z_0)>0$ at some point $z_0 \in \mathbb{D}$, then $v(z)=u_1(z)-u_2(z)>0$  in a neighborhood of
$z_0$. Thus 
$$ \Delta v(z)=\Delta u_1(z)-\Delta u_2(z)=4 \, |h(z)|^2 \, \left( e^{2
    u_1(z)}-e^{2 u_2(z)} \right) \ge 0 \, $$
there, i.e., $v$ is subharmonic in this neighborhood.
If $v(z_0)=0$ for some point $z_0 \in \mathbb{D}$, then $v$ satisfies the submean inequality
$$ v(z_0)=0 \le \frac{1}{2 \pi} \int \limits_{0}^{2 \pi} v\left(z_0+r \, e^{ i t}
\right) \, dt$$ for all $r$ small enough. Hence $v$ is subharmonic in $\mathbb{D}$.

\medskip

As we have already noted above, we have
\begin{equation} \label{eq:hi1}
 u_2(z) \ge \log \frac{1}{1-|z|^2} 
\end{equation}
from the Ahlfors--Yau lemma. An inequality in the opposite direction follows
from Theorem \ref{thm:liouville} and the Schwarz--Pick lemma. 
Namely, Theorem \ref{thm:liouville} gives us holomorphic functions $f_j : \mathbb{D}
\to \mathbb{D}$ such that
$$
 u_j(z)=\log \left( \frac{|f_j'(z)|}{1-|f_j(z)|^2} \frac{1}{|h(z)|} \right)
\, , 
  \qquad j=1,2 \, , 
$$
so the Schwarz--Pick lemma implies
\begin{equation} \label{eq:hi2}
 u_j(z) \le \log \left( \frac{1}{1-|z|^2} \frac{1}{|h(z)|} \right)\, ,
\qquad j=1,2\, .  
\end{equation}
Combining the estimates (\ref{eq:hi1}) and (\ref{eq:hi2}), we see that
\begin{equation} \label{eq:hi3}
 0 \le v(z) =\max \left\{u_1(z),u_2(z) \right\}-u_2(z) 
\le \log \frac{1}{|h(z)|} \, , \qquad z \in \mathbb{D} \, . 
\end{equation}
This, however, forces $v(z) \equiv 0$, because
the integral means
$\int_{0}^{2 \pi} v\left( r e^{i t} \right)  dt$
are monotonically increasing for $r \in (0,1)$ and
$$ 0 \le \int \limits_{0}^{2 \pi} v\left( r e^{i t} \right) \, dt
\overset{(\ref{eq:hi3})}{\le}
- \lim \limits_{r \nearrow 1} 
\int \limits_{0}^{2 \pi} \log \left| h(r e^{i t}) \right| \, dt = 0  \, $$
since $h$ is a Blaschke product (see \cite[Ch.~II, Theorem 2.4]{Ga}). Thus $v(z)
\equiv 0$ and $u_1\equiv u_2$. \hfill{$\square$}
\hfill{$\blacksquare$}

\subsection{Proof of Theorem \ref{thm:liouville}} \label{sec:liouville}

Before turning to the proof of Theorem \ref{thm:liouville}, we wish to
point out explicitly that
every $C^2$--solution $u : \Omega \to \mathbb{R}$ to $\Delta u=4 \, |h(z)|^2 \, e^{2 u}$
in $\Omega$ is in fact {\it real analytic} there. This follows from standard
elliptic regularity results for the Poisson equation and Bernstein's positive solution 
to Hilbert's Problem 19:
A priori the right--hand side of $\Delta u=4 \, |h(z)|^2 \, e^{2 u}$ is in
$C^2(\Omega)$, so $u$ certainly belongs to $C^3(\Omega)$ ({\it cf}.~\cite[Theorem
6.19]{GT97}). By Bernstein's analyticity theorem, $u$ is actually
real analytic in $\Omega$. This additional information will be a crucial
ingredient in the proof of Theorem \ref{thm:liouville}, which we now turn towards.

\medskip
It is convenient to make use of the $\partial/\partial z$-- and
$\partial/\partial\bar{z}$--operators,
$$ \frac{\partial}{\partial z}=\frac{1}{2} \left( \frac{\partial}{\partial x}-i \, 
  \frac{\partial}{\partial y} \right) \, , \qquad
 \frac{\partial}{\partial \bar{z}}=\frac{1}{2} \left( \frac{\partial}{\partial
     x}+i \, 
  \frac{\partial}{\partial y} \right) \, , \qquad z=x+i\, y\, . $$
In particular,
\begin{equation} \label{eq:laplace}
\Delta=\frac{\partial^2}{\partial x^2}+\frac{\partial^2}{\partial y^2}=4 \, 
\frac{\partial^2}{\partial z \, \partial \bar{z}} \, .
\end{equation}

We start off with the following simple, but important observation.
\begin{lemma} \label{lem:liou1}
Let $\Omega \subseteq \mathbb{C}$ be an open set, $h : \Omega \to \mathbb{C}$ 
a holomorphic function, $h \not\equiv 0$, and $u : \Omega \to \mathbb{R}$  a $C^2$--solution to
$\Delta u=4 \, |h(z)|^2 \, e^{2 u}$
in $\Omega$. Then the function
$$ B_u(z):=\frac{\partial^2 u}{\partial z^2}(z)-\left( \frac{\partial
    u}{\partial z}(z) \right)^2-\frac{h'(z)}{h(z)} \cdot \frac{\partial
  u}{\partial z}(z)$$
is holomorphic in $\Omega$ with the exception of possible simple poles at the
zeros of $h$.
\end{lemma}

{\bf Proof.}
Since we know that  $u$ is  of class $C^{\infty}$, we may certainly
differentiate the PDE $\Delta u=4 \, |h(z)|^2 \, e^{2 u}$ once with respect to $z$,
and  using (\ref{eq:laplace}), we get 
$$ \frac{\partial^3 u}{\partial z^2 \partial \bar{z}}=\frac{\partial^2
  u}{\partial z \partial \bar{z}} \frac{h'}{h}+\frac{\partial}{\partial
  \bar{z}} \left[ \left( \frac{\partial u}{\partial z} \right)^2 \right]  \, ,$$
so $\partial B_u/\partial \bar{z}=0$, which means
$B_u$ is meromorphic in $\Omega$.
\hfill{$\blacksquare$}

\medskip

In other words, $\frac{\partial u}{\partial z}$ is a formal (non holomorphic) solution of the
Riccati equation
$$ w'(z)-w(z)^2-\frac{h'(z)}{h(z)} \, w(z)= B_u(z) \, . $$
Following Laine \cite[p.~165]{Lai93} we transfer this Riccati equation via
$$ w(z)=v(z)-\frac{1}{2} \, \frac{h'(z)}{h(z)}$$
to normal form:
\begin{equation} \label{eq:riccati}
 v'-v^2= A_u(z) \, , 
\end{equation}
where
$$ A_u(z):=B_u(z)+ \frac{1}{2} \left( \frac{h'(z)}{h(z)} \right)'-\frac{1}{4}
\left( \frac{h'(z)}{h(z)} \right)^2 \, .$$
This function $A_u$ is holomorphic at every point $z_0 \in \Omega$ except when
$h(z_0)=0$. If $h$ has a zero of order $n \in \mathbb{N}$ at $z=z_0$, then
 $A_u$ has a pole of order $2$ there  with Laurent
expansion
\begin{equation} \label{eq:laurent1}
A_u(z)=\frac{b_0}{(z-z_0)^2}+\frac{b_1}{z-z_0}+b_2+ \cdots\, ,
\end{equation}
where
\begin{equation} \label{eq:laurent2}
 b_0=\frac{1-(n+1)^2}{4}\, . 
\end{equation}

Our next goal is to show that every local solution of the Riccati equation
(\ref{eq:riccati}) admits a
meromorphic extension to the whole of $\Omega$ provided $\Omega$ is a simply
connected domain.
 According to Laine
\cite[Theorem 9.1.7]{Lai93}
this is the case if and only if the first $n+2$ coefficients $b_0, \ldots, b_{n+1}$
in the Laurent expansion of $A_u$ satisfy a certain very complicated nonlinear relation
(see formula (6.18) in \cite{Lai93}), which appears to be difficult to verify
directly. We therefore choose a different path to establish:

\begin{lemma} \label{lem:liou2}
Let $\Omega \subseteq \mathbb{C}$ be a simply connected domain, $h : \Omega \to \mathbb{C}$ 
a holomorphic function, $h\not\equiv 0$, and $u : \Omega \to \mathbb{R}$  a $C^2$--solution to
$\Delta u=4 \, |h(z)|^2 \, e^{2 u}$
in $\Omega$. Then every local meromorphic solution to the Riccati equation
$v'-v^2=A_u(z)$
admits a meromorphic continuation to all of $\Omega$.
\end{lemma}

\begin{remark}
If $h(z) \equiv 1$ then Lemma \ref{lem:liou2}  reduces to the elementary
fact that if $A_u(z)$ is {\it holomorphic}, every solution to the Riccati
equation $v'-v^2=A_u(z)$ is meromorphic. Thus  Lemma \ref{lem:liou2}
is the essential step from Liouville's classical theorem for $\Delta
u=4 \, e^{2 u}$ to the more general Theorem \ref{thm:liouville}.
\end{remark}

{\bf Proof.}
Let $z_0 \in \Omega$ be a pole of $A_u$, i.e., 
a zero of order $n$ of $h$, so $A_u$ has an
expansion of the form (\ref{eq:laurent1})--(\ref{eq:laurent2}) at $z_0$.
We need only show that every local meromorphic solution to the Riccati equation
$v'-v^2=A_u(z)$ admits a meromorphic continuation to a neighborhood of $z_0$.
In order to simplify notation we take without loss of generality $z_0=0$.

\medskip

(i) We first show that $v'-v^2=A_u(z)$ has at least one meromorphic solution
in a neighborhood of the origin with residue $-(n+2)/2$. We substitute
$$ v_1(z)=\frac{- \frac{n+2}{2}+\omega(z)}{z} \,  ,$$
and find after some manipulation that $v_1'-v_1^2=A_u(z)$ if and only if
$$z\,  \omega'(z)=-(n+1) \, \omega(z)+\omega(z)^2+A_u(z) \, z^2
-\frac{1-(n+1)^2}{4} \, . $$
This is a Briot--Bouquet differential equation, which has a unique holomorphic
solution $\omega(z)$ in a neighborhood of $0$ such that $\omega(0)=0$, see
for instance \cite[Theorem 11.1.1]{Hil76}, because $-(n+1)$ is not a positive
integer. Therefore, the Riccati equation $v'-v^2=A_u(z)$ has at least one 
meromorphic solution $v_1$ with a simple pole at $0$ and such that the
residue of $2 v_1$ at $z_0$ is an integer.

\medskip

(ii) By Lemma 9.1.4 in \cite{Lai93}, we conclude that {\it every}
local meromorphic solution to the Riccati equation $v'-v^2=A_u(z)$
admits a meromorphic continuation to $0$ provided that we can exhibit a
second meromorphic solution in a neighborhood of $0$.

\medskip

(iii) Recall that 
$$\frac{\partial u}{\partial z}(z)+\frac{1}{2} \frac{h'(z)}{h(z)}$$
is a ``formal'' solution to $v'-v^2=A_u(z)$ with ``residue'' $n/2$ at $z=0$.
We extract an actual (meromorphic) solution $v_2$ from this formal solution as follows.
As $h$ has a zero of order $n$ at $z=0$, we have $h(z)=z^n
h_1(z)$ for a  function $h_1$ holomorphic at $z=0$ with
$h_1(0)\not=0$. A quick calculation shows that the function
$$ \nu(z):=u(z)+\log |h_1(z)| \, $$
satisfies 
\begin{equation} \label{Ahphi}
 A_u(z)=\frac{1-(n+1)^2}{4 z^2}+\frac{\partial^2 \nu}{\partial z^2} (z)-\frac{n}{z}\, 
\frac{\partial \nu}{\partial z}(z)-\left( \frac{\partial \nu}{\partial
    z}(z)\right)^2 \, .
\end{equation}
Recall that $u$ and hence also  $\nu$ is a {\it real analytic} function.
 This allows us to expand $\nu$ in a power series in $z$ and $\bar{z}$
in a neighborhood $U$ of $z=0$. We thus obtain for $ z, \bar{z} \in U$
$$ \nu(z,\bar{z})=\sum \limits_{k=0}^{\infty} \left(
\sum \limits_{j=0}^{\infty} a_{jk} z^j \right) \bar{z}^k=\sum
\limits_{j=0}^{\infty} a_{j0} z^j+\sum \limits_{k=1}^{\infty} \left(
\sum \limits_{j=0}^{\infty} a_{jk} z^j \right) \bar{z}^k\, ,$$
that is $$\nu(z, \overline{z})=g(z)+\Lambda(z,\bar{z})\, ,$$
if we set
$$ g(z)=\sum
\limits_{j=0}^{\infty} a_{j0} z^j \qquad \text{ and } \qquad \Lambda(z,\bar{z})=
\sum \limits_{k=1}^{\infty} \left(
\sum \limits_{j=0}^{\infty} a_{jk} z^j \right) \bar{z}^k.$$

%\begin{equation} \label{reprphi}
%\phi(z,\bar{z})=h(z)+\Lambda(z,\bar{z}).
%\end{equation}
Clearly,  $g(z)$ is holomorphic and $\Lambda(z,\bar{z})$ is real
analytic in $U$.

\medskip
Replacing $\nu$ by $g+\Lambda$ in (\ref{Ahphi}) yields
\begin{equation} \label{eq:q}
 A_u(z)=\frac{1-(n+1)^2}{4 z^2}- \frac{n}{z} \,
g'(z)-g'(z)^2+g''(z)+H(z,\bar{z})
\end{equation}
with
$$ H(z,\bar{z})=\Lambda_{zz}(z,\bar{z})-\frac{n}{z}\, 
\Lambda_z(z,\bar{z})-2 \, g'(z) \, \Lambda_z(z,\bar{z})-\left(
  \Lambda_{z}(z,\bar{z}) \right)^2.$$
Since $z H(z,\bar{z})$ is real analytic in $U$, $H(z,\bar{z})$ can be written
 as
$$ H(z,\bar{z})=\frac{1}{z} \cdot \left( \sum \limits_{k=1}^{\infty}
\left( \sum \limits_{j=0}^{\infty} b_{jk} z^{j} \right) \bar{z}^k
\right) \, . $$
Now, as $A_u$ has a pole of order $2$ at $z=0$, identity (\ref{eq:q}) shows that
$z^2 H(z,\bar{z})$ is 
{\it holomorphic} there:
$$ 0=\left( z^2 H(z,\bar{z}) \right)_{\bar{z}}=
\sum \limits_{k=1}^{\infty} \left( \sum \limits_{j=0}^{\infty} b_{jk} z^{j+1}
\right) k \bar{z}^{k-1}\, .$$
This implies $b_{jk}=0$ for all $j \in \mathbb{N}_0$ and $k \in \mathbb{N}$,
and consequently $H(z,\bar{z}) \equiv 0$.

\medskip
We thus obtain
$$ A_u(z)=\frac{1-(n+1)^2}{4 z^2}- \frac{n}{z} \,
g'(z)-g'(z)^2+g''(z) \quad  \text{for} \, \, z \in U \, , $$ 
and a  glance at the right hand side of this equation shows
$$ A_u(z)=\left( \frac{n}{2 z}+g'(z) \right)'-\left( \frac{n}{2
    z}+g'(z) \right)^2.$$
Therefore the function
$$ v_2(z)=\frac{n}{2z}+g'(z)$$
is a meromorphic solution of the Riccati differential equation
$v'-v^2=A_u(z)$ in $U$. 
\hfill{$\blacksquare$}

\medskip

{\bf Proof of Theorem \ref{thm:liouville}.}
Let $u$ be a $C^2$--solution to $\Delta u=4 |h(z)|^2 \, e^{2 u}$ in $\Omega$
and $z_0$ {\it not} a zero of $h(z)$. Then
$$ u_1(z):=u(z)+\log|h(z)|$$
is a $C^2$--solution to Liouville's equation $\Delta u_1=4 e^{2 u_1}$ in a
neighborhood $V$ of $z_0$. By Liouville's theorem, we get
$$ u_1(z)=\log \frac{|f'(z)|}{1-|f(z)|^2}$$
for a holomorphic function $f : V \to \mathbb{D}$, i.e.~(\ref{eq:liouville}) holds
for $z \in V$. A straightforward calculation now shows that the Schwarzian of $f$ 
$$ S_f(z):=\left( \frac{f''(z)}{f'(z)} \right)'-\frac{1}{2} \left(
  \frac{f''(z)}{f'(z)} \right)^2$$
satisfies
\begin{equation} \label{eq:schwarzian}
S_f(z)=2 A_u(z)
\end{equation}
in $V \subseteq \Omega$. By Lemma \ref{lem:liou2}, the associated Riccati equation
$v'-v^2=A_u(z)$ has {\it only} meromorphic solutions in $\Omega$, so
the Schwarzian differential equation (\ref{eq:schwarzian}) {\it only} has meromorphic
solutions in $\Omega$ as well. This follows from
Theorem 9.1.7 and Corollary 6.8 in \cite{Lai93}. In particular, $f : V \to \mathbb{D}$
has a meromorphic continuation to $\Omega$ which we continue to call $f$. 
We finally note that $|f(z)|<1$ in all of $\Omega$. In fact,
let $\Omega'$ be the component of $\{z \in \Omega \, : |f(z)|<1\}$ which contains
$V$. Then $\Omega'$ is clearly open in $\Omega$, but also
 closed in $\Omega$. This is immediate from the fact that 
(\ref{eq:liouville}) holds for all $z \in \Omega'$. 

\medskip

In order to prove the uniqueness statement, let $f, g : \Omega \to \mathbb{D}$ be two
holomorphic maps such that
\begin{equation} \label{eq:uniq}
u(z)= \log \left( \frac{|f'(z)|}{1-|f(z)|^2} \frac{1}{|h(z)|} \right)
=\log \left( \frac{|g'(z)|}{1-|g(z)|^2} \frac{1}{|h(z)|} \right) \, . 
\end{equation}
In particular, $f$ and $g$ are nonconstant and
as above  ${\cal S}_f(z) ={\cal S}_g(z)$ in $\Omega$, so
$f=T \circ g$ for some M\"obius transformation $T$. Thus (\ref{eq:uniq}) shows
$$ \frac{|T'(w)|}{1-|T(w)|^2}=\frac{1}{1-|w|^2}$$
first for all points $w$ in the open set $g(\Omega) \subseteq \mathbb{D}$
 and then clearly for every $w \in \mathbb{D}$. Hence $T(\mathbb{D}) \subseteq \mathbb{D}$ and
the Schwarz--Pick lemma implies that $T$ is a conformal disk automorphism.
\hfill{$\square$}
\hfill{$\blacksquare$}

\begin{remark}
The above proof of Theorem \ref{thm:liouville} uses Liouville's theorem (i.e.,
the special case $h(z) \equiv 1$).  This can be avoided by showing directly 
as in \cite{Min} or \cite{Nit57} that the
solution $f$ to the initial value problem
$$ {\cal S}_f(z)=2 A_u(z) \, , \qquad f(z_0)=0 \, , \quad f'(z_0)=e^{u_1(z_0)} \, ,
\quad f''(z_0)=2 e^{u_1(z_0)} \frac{\partial u_1}{\partial z}(z_0) \, ,$$
which is meromorphic in all of $\Omega$ by Lemma \ref{lem:liou2}, 
fulfills (\ref{eq:liouville}) in a neighborhood of $z_0$.
In particular, the equation (\ref{eq:liouville}) can be solved for $f$ {\it constructively}.
\end{remark}

\subsection{Proof of Theorem \ref{thm:4} and  Theorem \ref{thm:5}}

{\bf Proof of Theorem \ref{thm:4}.}
Let $p$ be a polynomial with  zeros $z_j$ (counted with multiplicities).
In view of Lemma \ref{lem:help} (d)
there exists a uniquely determined  
conformal Riemannian metric $\lambda(z)\, |dz|$ in $\Omega$ with curvature
$-4 \, |p(z)|^2$ and boundary values $ \phi(\xi)/ |p(\xi)|$. From Theorem 
\ref{thm:liouville} we deduce
$$\lambda(z)=\frac{|f'(z)|}{1-|f(z)|^2} \frac{1}{|p(z)|} \, , \qquad  z \in \Omega,$$
for some holomorphic function $f : \Omega \to \mathbb{D}$. Thus $\{ z_j \}$ is the set
of critical points of  $f$, and the boundary condition 
(\ref{eq:bdd01}) is fulfilled, because
$$ \lim \limits_{z \to \xi} \frac{|f'(z)|}{1-|f(z)|^2}
=\lim \limits_{z \to \xi} \lambda(z) \, |p(z)|=\phi(\xi) \, .
$$
If $g$ is another holomorphic function $g : \Omega \to \mathbb{D}$ with the properties
stated in Theorem \ref{thm:4}, then
$$\tilde{\lambda}(z):=\frac{|g'(z)|}{1-|g(z)|^2} \frac{1}{|p(z)|}$$
is the density of a regular conformal Riemannian metric in $\Omega$ of curvature
$-4 \, |p(z)|^2 $ in $\Omega$ and boundary values $\phi/|p|$. 
From the uniqueness statement of
Lemma \ref{lem:help} (d) we infer $\tilde{\lambda}=\lambda$ in $\Omega$, that is
$$ \frac{|g'(z)|}{1-|g(z)|^2}\frac{1}{|p(z)|}=\frac{|f'(z)|}{1-|f(z)|^2}
\frac{1}{|p(z)|} \quad \text{for} \,
\, z \in \Omega\, . $$
Hence, applying Theorem \ref{thm:liouville}, we see that $g=T \circ f$ for
some conformal automorphism $T$ of $\mathbb{D}$.  \hfill{$\square$}
\hfill{$\blacksquare$}

{\bf Proof of Theorem \ref{thm:5}.}
We only prove the existence part.
Let $B$ be a Blaschke product with zeros $z_j$ (counted with multiplicities).
Note, $$\lim_{r \to 1-} |B(r \xi)|=1  \text{ for a.e. } \xi \in \partial
\mathbb{D}\,. $$
Next, let $v$ be the harmonic function in $\mathbb{D}$ with boundary values 
$\log \phi \in L^{\infty}(\partial \mathbb{D})$, so $|v(z)| \le M$ in $\mathbb{D}$ for some
constant $M>0$ and 
$$ \angle \lim \limits_{z \to \xi} v(z)=\log \phi(\xi) \qquad \text{for a.e. }
\xi \in \partial \mathbb{D}\, .$$
Then $v(z)=\Re \log g(z)$ for some nonvanishing holomorphic function $g : \mathbb{D}
\to \mathbb{C}$.

\medskip

By Lemma \ref{lem:help} (d) there exists a unique conformal Riemannian metric
$\mu(z) \, |dz|$ with curvature $-4 \, |B(z)\, g(z)|^2 $ in $\mathbb{D}$ and $\mu(\xi)=1$
for $\xi \in \partial \mathbb{D}$. Thus we can apply Theorem \ref{thm:liouville}
and get a holomorphic function $f : \mathbb{D} \to \mathbb{D}$ with critical points $z_j$ (counted with
multiplicities) and no others and
$$ \mu(z)=\frac{|f'(z)|}{1-|f(z)|^2} \frac{1}{|B(z)| \, |g(z)|} \quad \text{ for } \, \,
z \in \mathbb{D}.$$
By construction,
$$ \angle \lim \limits_{z \to \xi} \frac{|f'(z)|}{1-|f(z)|^2}=
\angle \lim \limits_{z \to \xi} \mu(z) \, |B(z)| \, |g(z)| =
\phi(\xi)
  \quad \text{ for a.e. } \xi \in  \partial \mathbb{D}\, $$
and
\[
 \sup \limits_{z \in \mathbb{D}} \frac{|f'(z)|}{1-|f(z)|^2} <\infty
 \, .
\]
\hfill{$\blacksquare$}

dakraus@mathematik.uni-wuerzburg.de\\
roth@mathematik.uni-wuerzburg.de
\end{document}